\documentclass[12pt,reqno]{amsart}
\usepackage{amscd,amssymb,amsthm,euscript}

\newtheorem{thm}{Theorem}[section]
\newtheorem{prop}[thm]{Proposition}
\newtheorem{lem}[thm]{Lemma}
\newtheorem{cor}[thm]{Corollary}
\newtheorem*{theorem}{Theorem}
\newtheorem{defn}[thm]{Definition}
\theoremstyle{remark}

\newtheorem{ex}[thm]{Example}
\newtheorem*{rem}{Remark}

\newcounter{ictr}
\newenvironment{ilist}{\begin{list}
                         {\textup{(\roman{ictr})}}
                         {\usecounter{ictr}
                          \setlength{\leftmargin}{0.6truein}
                          \setlength{\itemsep}{0.0truein}
                          \setlength{\labelwidth}{0.3truein}}}
                      {\end{list}}
\newcommand{\Hawaii}{Hawai\kern.05em`\kern.05em\relax i}
\newcommand{\Manoa}{M\=anoa}
\DeclareMathOperator{\isom}{\ensuremath{\mathop{Isom}}}
\DeclareMathOperator{\op}{\ensuremath{\mathop{Op}}}
\DeclareMathOperator{\LE}{\ensuremath{\mathop{Lip^{ls}}}}

\DeclareMathOperator{\LEG}{\ensuremath{\mathop{Lip^{ls}_\Gamma}}}
\DeclareMathOperator{\card}{\ensuremath{\mathop{card}}}
\renewcommand{\epsilon}{\varepsilon}
\newcommand{\eps}{\varepsilon}
\renewcommand{\phi}{\varphi}

\newcommand{\ga}{\Gamma}
\newcommand{\crg}{C_r^*(\ga)}
\newcommand{\F}{\mathbb{F}}
\newcommand{\N}{\mathbb{N}}
\newcommand{\R}{\mathbb{R}}
\newcommand{\Z}{\mathbb{Z}}
\newcommand{\C}{\mathbb{C}}
\newcommand{\h}{\mathcal{H}}
\newcommand{\cc}{\colon}
\newcommand{\cont}{\subseteq}
\newcommand{\x}{\times}

\begin{document}

\title[Exactness and Uniform Embeddability]%
    {Exactness and Uniform Embeddability \\of Discrete Groups}
\author{Erik Guentner}
\address{Mathematics Department, University of \Hawaii, \Manoa, 2565 McCarthy
   Mall, Honolulu, HI 96802}
\email{erik@math.hawaii.edu}
\author{Jerome Kaminker}
\address{Department of Mathematical Sciences, IUPUI, 402 N.~Blackford
   Street, Indianapolis, IN 46202-3216} 
\email{kaminker@math.iupui.edu}

\thanks{The first author was supported in part by an MSRI Postdoctoral
  Fellowship and NSF Grant DMS-0071402.  The second author was supported
  in part by NSF Grant DMS-0071435}  

 
\begin{abstract}
  We define a numerical quasi-isometry invariant, $R(\Gamma)$, of a
  finitely generated group 
  $\Gamma$, whose values parametrize the difference between
  $\Gamma$ being  uniformly embeddable in a Hilbert space
  and $C^{*}_{r}(\Gamma)$ being exact.
\end{abstract}
\maketitle

\section{Introduction}

In his study of large scale properties of finitely generated groups,
M. Gromov introduced the notion of uniform embeddability.
\cite{gromov93}.  Recall that a {\it uniform embedding\/} of one
metric space $(X,d_{X})$ into another $(Y,d_{Y})$ is a function $f \cc
X \to Y$ for which there exist non-decreasing functions $\rho_{\pm}
\cc [0,\infty) \to \R$ such that $\lim_{r\to\infty}\rho_\pm(r)=+\infty$
and such that for all
$x$, $y\in X$
\begin{equation}
\label{eqn:ue}
\rho_{-}(d_{X}(x,y)) \leq d_{Y}(f(x),f(y)) \leq \rho_{+}(d_{X}(x,y)).  
\end{equation}
The condition 
$\lim_{r\to\infty}\rho_\pm(r)=+\infty$ is summarized by saying that the 
$\rho_\pm$ are {\it proper}.
In an appendix we collect several known facts about the relation
between uniform embeddings and other notions from coarse geometry.

Gromov raised the question of whether a finitely
generated group
that is uniformly embeddable in a Hilbert space (when viewed as a metric
space with a word length metric) 
satisfies the Novikov Conjecture \cite{NCITR}.  This was answered
affirmatively by Yu:

\begin{theorem}[\cite{yu00,skandalis-tu-yu02}]
Let $\Gamma$ be a finitely generated group, equipped with a word length
metric.  If $\ga$ is uniformly
embeddable in Hilbert space then $\ga$ satisfies both the Novikov
Conjecture and the Coarse Baum-Connes Conjecture. 
\end{theorem}

Recently, Gromov has proved the existence of a countable discrete group which is not
uniformly embeddable  in a Hilbert space \cite{gromov:random}.
On the other hand, it  has been observed that this group does
satisfy the Novikov Conjecture, although it is not known whether it
satisfies the Coarse Baum-Connes Conjecture
\cite{higson-lafforgue-skandalis02}.

From the analytic side, E. Kirchberg and S. Wassermann extensively
studied the notion of exactness of a countable discrete group
\cite{kirchberg93,kirchberg-wassermann99b}.  Recall that $\ga$ is {\it
  exact\/} if $\crg$ is an exact $C^{*}$-algebra, that is, if taking
minimal tensor product with $\crg$ on each of the terms in a short
exact sequence of $C^*$-algebras preserves the exactness of the
sequence.

Uniform embeddability is a geometric property of a group, while
exactness is more closely related to harmonic analysis.
It is interesting that there is a relation between these notions.
Indeed, the connection between these types of properties is related to the
Baum-Connes Conjecture.  Recently, it was shown that exactness of a
countable discrete group implies its uniform embeddability in a
Hilbert space.

\begin{theorem}[\cite{guentner-kaminker02,guentner-kaminker02a,ozawa00}]
  Let $\ga$ be a finitely generated discrete group.  If
  $C_{r}^{*}(\Gamma)$ is an exact $C^{*}$-algebra, then $\Gamma$, viewed
  as a metric space with a word length metric, is uniformly embeddable
  in a Hilbert space.
\end{theorem}

One may ask to what extent the converse of this result holds.  The
question of whether a uniformly embeddable group is exact has been
studied from various perspectives (see
\cite{dadarlat-guentner-amalgam,guentner-kaminker:cbgrowth}).  In the
present paper we introduce a numerical invariant $R(\Gamma)$ of a
finitely generated discrete group $\ga$ which can be viewed as
parameterizing the difference between the group being exact and being
uniformly embeddable in a Hilbert space.

\section{The definition of $R(\Gamma)$ }

Although our primary interest is in uniform embeddings into Hilbert
space, we will formulate the basic definitions in the context of general 
metric spaces.  
Recall that a function $f:X \to Y$ is {\it large-scale Lipschitz\/} if
there exist $C>0$ and $D\geq 0$ such that
\begin{eqnarray}
\label{eq:lslip}
  d_Y(f(x),f(y)) \leq C d_X(x,y) + D.
\end{eqnarray}
The following example shows that a uniform embedding of a discrete metric
space is not necessarily Lipschitz, or even large-scale Lipschitz.

\begin{ex}
\label{ex0}
  Let $X=\{\; (n,1/n),\,(n,0) : \text{$n=1,2,\dots$} \;\} \subset
  \R^2$ with the induced metric.  Define $f:X\to\R^2$ by $f(n,1/n) =
  (n,1)$ and $f(n,0)=(n,0)$.  Then $f$ is a both a uniform embedding
  and a large-scale Lipschitz map but  it is not a Lipschitz
  map.  Let $Y=\{\, n^2 : \text{$n\in\N$} \;\} \subset \R$ with the
  induced metric.  Define $g:Y\to\R$ by $g(y)=y^2$.  Then $g$ is a
  uniform embedding but it is not large-scale Lipschitz.
\end{ex}

Let $\LE(X,Y)$ denote the set of large-scale Lipschitz maps from $X$
to $Y$.  Following Gromov \cite{gromov93}, define
the {\it compression\/} $\rho_f$ of $f \in\LE(X,Y)$ by
\begin{equation}
\label{eqn:c}
  \rho_{f}(r) = \inf_{d_X(x,y) \geq r} d_Y(f(x),f(y)) 
\end{equation}
The compression function $\rho_{f}$ is a non-decreasing, non-negative
real-valued function satisfying the first inequality in (\ref{eqn:ue}),
and has the property that if $\rho_{-}$ is a another
such function then
$\rho_- \leq \rho_f$.  Consequently, $f$ is a uniform embedding if and
only if $\rho_{f}$ is proper.  Always assuming that the metric on $X$
is unbounded, we define a real-valued invariant of $X$
 as follows.

\begin{defn}
\label{eqn:rgamma}
Let $X$ be a metric space with an unbounded metric.
\begin{ilist}
\item  The {asymptotic compression} $R_f$
of a large-scale Lipschitz map $f\in\LE(X,Y)$ is
\begin{equation}
\label{eqn:asc}
R_{f}= \liminf_{r \to \infty}
    \displaystyle{\frac{\log\rho_{f}^*(r)}{\log r}},
\end{equation}
where $\rho_f^*(r) = \max\{\, \rho_f(r),1 \,\}$.
\item The {compression} of X in Y is 
\begin{equation*}
R(X,Y) = \sup\{\, R_{f} : \text{$f \in \LE(X,Y)$} \,\}.
\end{equation*}
\item If Y is a Hilbert space, then
the {Hilbert space compression\/} of X is
\begin{equation*} 
R(X) = R(X,\h).
\end{equation*}
\end{ilist}
\end{defn}

\begin{rem}
  The distinction between $\rho_f^*$ and $\rho_f$ in (\ref{eqn:asc}) is
  not essential, but 
  is meant to eliminate pathology; when
  $\rho_f$ is unbounded, the definition (\ref{eqn:asc}) is unchanged if
  we replace $\rho_f^*$ by $\rho_f$.   Also, observe that $R_f\geq 0$. 
\end{rem}

\begin{prop}
\label{ex1}
The  compression of  $X$ in $Y$ satisfies $R(X,Y)\leq 1$.  Indeed, the
asymptotic compression of a large-scale Lipschitz map $f$ satisfies
$R_f\leq 1$.  
\end{prop}
\begin{proof}
Let $f \in \LE(X,Y)$ and let
$C>0$ and $D\geq 0$ be constants as in (\ref{eq:lslip}) supplied by
fact that 
$f$ is large-scale Lipschitz.  Since $X$ is unbounded, there exist
sequences $x_n$ and $y_n\in X$ such that $r_n=d_X(x_n,y_n)\to\infty$.
For these $r_n$ we have $\rho_{f}(r_n) = 
\inf_{d(x,y) \geq r_n} d_{Y}(f(x), f(y)) \leq Cr_n + D$ and, for all
sufficiently large $n$, $\rho_f^*(r_n)\leq Cr_n +D$.  Hence
\begin{equation*}
R_{f} = \liminf_{r \to \infty}\displaystyle{\frac{\log\rho^*_{f}(r)}{\log r}} 
   \leq \liminf_{n \to \infty}\displaystyle{\frac{\log(Cr_n+D)}{\log r_n}} 
          = 1.  \quad\qed
\end{equation*}
\renewcommand{\qed}{}
\end{proof}

\begin{prop}
\label{ex2}
If a metric space $X$ admits an isometric embedding into
a metric space $Y$, then $R(X,Y) = 1$.
\end{prop}
\begin{proof}
If $f:X\to Y$ is an isometry, we have $\rho_{f}(r) \geq r$, hence $R_{f} =1$.
Thus, $R(X,Y) \geq 1$.  By the previous proposition $R(X,Y) = 1$. 
\end{proof}

In fact, the same conclusion will follow from the
existence of a quasi-isometric embedding of $X$ into a Hilbert space (see
Theorem~\ref{thm:qii}).

\begin{prop}
If a metric space $X$ admits an isometric embedding into the Banach
space $l^1(\N)$ then $R(X) \geq 1/2$.
\end{prop}
\begin{proof}

  Let $f: X \to l^1(\N)$ be an isometric embedding.  Define a
  function $g: \R\to L^2(\R)$ by mapping $x\geq 0$ to the
  characteristic function of $[0,x]$ and $x<0$ to the characteristic
  function of $[x,0]$.  Note that $\| g({x})-g({y})
  \|^2_{l^2(\R)} = |{x} -{y}|$.  For $\mathbf{x}=(x_1,x_2,\dots)\in
  l^1(\N)$ set
  $h(\mathbf{x}) =g(x_1)\oplus g(x_2) \oplus \cdots\in L^2(\R)\oplus
  L^2(\R)\oplus\cdots$.   Then it is
  easily checked that $h\circ f$ is a uniform embedding
  with $\rho^{*}_{f}(r) \geq \sqrt{r}$.  Hence, $R(X) \geq 1/2$.
\end{proof}

 The same conclusion would follow from the existence of a
quasi-isometric embedding of $X$ into $l^1(\N)$.

We record here some results which will be proved later in the paper.

\begin{ex} 
  The Hilbert space compression of the metric space obtained from a
  sequence of expander graphs is zero, as is true of any metric space
  that is not uniformly embeddable in a Hilbert space (see
  Proposition~\ref{prop:ue}).
\end{ex}

\begin{ex}
\label{ex4}
According to Example~\ref{ex2} we have $R(\Z)=1$.  Further,
$R(\Z^n)=1$, for all $n\in\N$ (see Proposition~\ref{prod}).
\end{ex}

\begin{ex}
  Let $\F_2$ be the free group on two generators.  We have $R(\F_2)=1$
  (see Proposition~\ref{freegp}).
\end{ex}

Our primary interest is when our metric space is a finitely generated
discrete group $\ga$, equipped with the left invariant metric induced by
the {\it word length\/} function associated to a finite, symmetric
generating set.

The value $R(\Gamma)$ will be shown to be independent of the particular
generating set chosen.  
Equipped with a metric in this manner, $\ga$ is a geodesic space.  In fact it
will be important to have a result that holds even for a quasi-geodesic
space.

Recall that a discrete metric space $X$ is a {\it quasi-geodesic
  space\/} if there exist $\delta>0$ and $\lambda \geq 1$ such that for
all $x$ and $y\in X$ there exists a sequence $x=x_0, x_1,\dots,x_n =y$
of elements of $X$ such that
\begin{align}
\label{eqn:qg}
\begin{split}
\sum_{1}^{n} d_X(x_{i-1},x_{i}) &\leq \lambda \,d_X(x,y),   \\
d_X(x_{i-1},x_{i}) &\leq \delta, \quad\text{for all $1\leq i \leq n$.}      
\end{split}
\end{align}
Although a uniform embedding of a discrete metric space is
not necessarily large-scale Lipschitz (Example~\ref{ex0}), a uniform
embedding of a quasi-geodesic metric space is.

\begin{prop}[\cite{gromov93}]
\label{qg}
  Let $X$ and $Y$ be metric spaces, and assume that $X$ is  
  quasi-geodesic.  Let $f : X \to Y $ be a uniform embedding %
  Then $f$ is large-scale Lipschitz.
\end{prop}
\begin{proof}
We will only use the existence of $\rho_+$ (a non-decreasing,
non-negative real-valued function satisfying the second inequality
in  (\ref{eqn:ue})).  Let $\lambda\geq 1$ and $\delta>0$ be the
constants supplied by the fact 
that $X$ is quasi-geodesic.  We will show that there exist
constants $C>0$, $D\geq 0$ such that 
\begin{eqnarray*}
    d_Y(f(x),f(y)) \leq C d_X(x,y)+D,\quad\text{for all $x$, $y\in X$}.
\end{eqnarray*}
Let $x$, $y\in X$ and let $x_0,\dots,x_n$ be a sequence of elements of
$X$ satisfying (\ref{eqn:qg}). 
Extract a subsequence $x_{i_0},\dots,x_{i_m}$ as follows: $i_0=0$
and, assuming $i_0,\dots,i_{j}$ are already defined,
\begin{equation*}
  i_{j+1} = \left\{ 
    \begin{minipage}{3.5truein}
         the smallest integer $k$ such that
         $d(x_{i_{j}},x_k)\geq \delta/2$, if such exists;
         if no such $k$ exists put $m=j$ and stop.
    \end{minipage}\right.
\end{equation*}
The subsequence has the following properties:

\begin{ilist}
  \item $x_{i_0}=x$, $d_X(x_{i_m},y)\leq \delta/2$, and
  \item $\delta/2 \leq d_X(x_{i_{j-1}},x_{i_{j}})\leq 3\delta/2$, 
               for $1\leq j \leq m$.
\end{ilist}

We have the following estimates:
\begin{equation}
  \begin{split}
  d_Y(f(x),f(y)) &\leq \sum_{j=1}^{m} d_Y(f(x_{i_{j-1}}),f(x_{i_j})) 
                   + d_Y(f(x_{i_m}),f(y))
             \leq m \rho_+(3\delta/2) + \rho_+(\delta/2), \\
  m\delta/2 &\leq \sum_{j=1}^m d_X(x_{i_{j-1}},x_{i_j}) \leq 
       \sum_{i=1}^n d_X(x_{i-1},x_i) \leq \lambda d_X(x,y).
  \end{split}
\end{equation}
From the second we conclude that $m\leq 2\delta^{-1} \lambda d_X(x,y)$
which, combined with the first, yields
\begin{eqnarray*}
  d_Y(f(x),f(y)) \leq  2\delta^{-1}\lambda \rho_+(3\delta/2) \, d_X(x,y) 
        + \rho_+(\delta/2). \qed
\end{eqnarray*}
\renewcommand{\qed}{}
\end{proof}

We next establish the quasi-isometry invariance of the Hilbert space
compression.  Recall that a function $\phi \cc X\to Y$ is a
{\it quasi-isometry \/} if there exist $C>0$ and $D\geq 0$ such that
for all $x$, $x'\in X$
\begin{eqnarray}
\label{eq:qidef}
  C^{-1} d_X(x,x')-D \leq d_Y(f(x),f(x'))\leq C d_X(x,x') +D,
\end{eqnarray}
The spaces $X$ and $Y$ are {\it quasi-isometric\/} if there exists a
quasi-isometry $\phi \cc X\to Y$ and $K>0$ such that $\phi$ has
$K$-dense range, meaning that every element of $Y$ is within distance
$K$ of an element in the image of $\phi$.

Equivalently, $X$ and $Y$ are quasi-isometric if there exist
quasi-isometries $\phi \cc X\to Y$ and $\psi \cc Y\to X$ and a $K>0$
such that 
\begin{equation*}
\begin{split}
  d(\psi\phi(x),x) &\leq K,\quad\text{for all $x\in X$} \\
  d(\phi\psi(y),y) &\leq K,\quad\text{for all $y\in Y$.}    
\end{split}
\end{equation*}
In this case, one calls $\phi$  a quasi-isometric equivalence.

\begin{prop}
\label{qi1}
 Let $\phi:X_{1}\to X_{2}$ be a quasi-isometry.  Then $R_{\phi} = 1$ 
\end{prop}
\begin{proof}
By Proposition ~\ref{ex1} we have $R_{\phi} \leq 1$.
If $C$ and $D$ are constants as in (\ref{eq:qidef}) supplied by the fact
that $\phi$ is a 
quasi-isometry we have
$C^{-1}r -D \leq \rho_{\phi}(r)$,
from which we conclude that $R_{\phi} \geq 1$.  
\end{proof}

\begin{prop}
\label{mult}
  Let $f \in \LE(X,Y)$ and $g \in \LE(Z,X)$.  Then 
  $f\circ g\in \LE(Z,Y)$ and $R_{f\circ g} \geq R_{f} R_{g}$.
\end{prop}
\begin{proof}
Let $f$ and $g$ be as in the statement.
Direct computation shows that $f\circ g$ is large-scale Lipschitz, and
further that 
  \begin{equation*}
    \rho_{f\circ g}(r) \geq \rho_{f}(\rho_{g}(r)),
\end{equation*}
and the same for $\rho^*$.
If the increasing function $\rho_{g}(r)$ is bounded then $R_g=0$ and
there is nothing to prove.  We therefore may assume that 
$\lim_{r\to \infty} \rho_{g}(r) = +\infty$.  From this and the previous
inequality we conclude that
\begin{equation*}
  \begin{split}
  R_{f\circ g} &= \liminf_{r\to\infty} 
        \displaystyle{
      \left\{\frac{\log \rho^*_{f\circ g}(r)}{\log r}\right\}} \\
    &\geq \liminf_{r\to\infty}
   \displaystyle{\left\{\frac{\log \rho^*_f(\rho_{g}(r))}{\log (\rho_{g}(r))}\right\}
           \left\{\frac{\log (\rho_{g}(r))}{\log r} \right\}} \\
    &\geq R_f R_{g}. \qquad \qed
\end{split}
\end{equation*} \renewcommand{\qed}{}
\end{proof}

\begin{thm}
\label{thm:qii}
Let $X_{1},X_{2}$ be metric spaces.  If there exists a quasi-isometry 
$\phi : X_{1} \to X_{2}$ then $R(X_{1},Y)\geq R(X_{2},Y)$, for every
metric space $Y$.
\end{thm}
\begin{proof}
Let $\phi:X_{1}\to X_{2}$ be a quasi-isometry, and let $Y$ be a metric
space.  If $f \in \LE (X_{2},Y)$ then $f\circ \phi \in \LE (X_{1},Y)$
and it follows from  
Propositions~\ref{mult} and \ref{qi1}  that 
$R_{f\circ \phi} \geq R_{f}R_{\phi} = R_{f}$.
Thus we have

\begin{equation*}
\begin{aligned}
   R(X_{1},Y) &= \sup \{\, R_{g} : g \in \LE(X_{1},Y) \,\} \geq 
              \sup \{\, R_{f\circ \phi} : f \in \LE(X_{2},Y) \,\} \\ 
       &\geq  \sup \{\, R_{f} : f \in \LE(X_{2},Y) \,\} = R(X_{2},Y). \quad\qed
\end{aligned}
\end{equation*}

\renewcommand{\qed}{}\end{proof}

\begin{cor}
\label{cor:qi-inv}
  If the metric spaces $X_{1}$ and $X_{2}$ are quasi-isometric then 
  $R(X_{1},Y)=R(X_{2},Y)$ for every metric space $Y$. \qed
\end{cor}

\begin{cor}
\label{cor:fgg}
  Let $\ga$ be a finitely generated discrete group.  The Hilbert space
  distortion 
  $R(\Gamma)$ is independent of the finite, symmetric generating set
  used to define the length function and metric on $\Gamma$.
\end{cor}
\begin{proof}
  Word length metrics associated to finite, symmetric generating
  sets are quasi-isometric; indeed, the identity provides the required
  $K$-dense quasi-isometry  \cite{delaharpe-ggt}.
\end{proof}

\section{Uniform embeddings and exactness}

In this section we will relate the Hilbert space compression of a
metric space $X$ to uniform embeddability, and,  in the case of
a finitely generated discrete group, to exactness.

\begin{prop}
\label{prop:ue}
Let $X$ be a metric space.  If the Hilbert space compression of $X$ is
nonzero then $X$ is uniformly embeddable in Hilbert space.
\end{prop}
\begin{proof}
Let $X$ be given with $R(X) > 0$.  From the
Definition~\ref{eqn:rgamma} of $R(X)$ we see that there exists
$\epsilon > 0$ and a large scale Lipschitz map $f \in \LE(X,\h)$ with
 asymptotic compression greater than $\epsilon$:
\begin{equation*}
  R_{f} = \liminf_{r \to\infty}
    \displaystyle{\frac{\log\rho^*_{f}(r)}{\log r}} > \epsilon.
\end{equation*}
In particular, for all sufficiently large $r$ , we have 
$\log \rho^*_{f}(r) \geq \frac{\epsilon}{2}\log r$, hence 
$\rho^*_{f}(r) \geq r^{\epsilon /2}$.  Consequently $\rho^*_{f}$, and
$\rho_f$, are proper
and $f$ is a uniform embedding.
\end{proof}

The main result of this section is the following.

\begin{thm}
\label{thm:exact}
  Let $\Gamma$ be a finitely generated discrete group.  If the Hilbert
  space compression of $\Gamma$ is greater than $1/2$ then $\Gamma$ is
  exact.
\end{thm}

The proof of the theorem
relies on the following characterization of exactness 
\cite{guentner-kaminker02,guentner-kaminker02a,ozawa00}:

\begin{prop}
\label{exact}
Let $\ga$ be a finitely generated discrete group, equipped with word
length and metric associated to a finite, symmetric set of
generators.  Then $\ga$ is exact if and only if there exists a
sequence  of  positive definite functions, $u_{n :}\ga\times\ga\to\R$, 

satisfying 

\begin{equation}
\label{convergence}
  \text{${\text for\ all\ \ } C>0$, $u_n\to 1$ uniformly on the strip 
         $\{\, (s,t) : d(s,t)\leq C \,\}$}
\end{equation}
and 
\begin{equation}
\label{support}
  \text{$for \ all\  n,\  there\  exist \ R > 0$,  such that $u_n(s,t)=0$ if
            $d(s,t)\geq R$}. \quad\qed
\end{equation}
\renewcommand{\qed}{}\end{prop}

\noindent  We refer to (\ref{convergence}) as the {\it convergence
  condition\/} and to (\ref{support}) as the {\it support condition\/};
a kernel $\ga\times\ga\to \R$ satisfying the support condition is of
{\it finite width}.

Under the assumption that $R(\Gamma) > 1/2$ we will construct a sequence of
positive definite kernels on $\ga\times\ga$ satisfying the convergence
and support conditions of the proposition.

Given a
complex-valued kernel $k:\Gamma \times \Gamma \to \C$, define an
operator $\op(k)$ by  convolution:
\begin{equation}
\label{eqn:convolve}
  \op(k)\xi(x) = \sum_{y\in Y} k(x,y) \xi(y), 
                      \quad\text{$\xi\in l^2(\Gamma)$.}
\end{equation}
We will need both of the following criteria for the boundedness of $Op(k)$ on
$l^{2}(\Gamma)$, (c.f. \cite{pedersen:analysis}). 

\begin{prop} Under either of the following conditions $Op(k)$ is a
  bounded operator.
\begin{ilist}
\item If $k$ is bounded and has finite width then  $\op(k)$ is bounded  
\item \textup{(}Schur Test\textup{)} Let 
$k$ be non-negative and real-valued  with the  property that
there exists $C>0$ such that
\begin{equation}
\label{eqn:schur}
\begin{split}
  \sum_{s\in \Gamma} k(s,t) &\leq C, \quad\text{for all $t\in \ga$} \\
  \sum_{t\in \Gamma} k(s,t) &\leq C, \quad\text{for all $s\in \ga$}.   
\end{split}
\end{equation}
Then  $\op(k)$ is bounded and
 $\| \op(k) \|\leq C$.
\end{ilist}   
\end{prop}

\begin{proof}[Proof of Theorem~\ref{thm:exact}]
Let $\Gamma$ be a finitely generated discrete group equipped with the
word length metric associated to a finite
symmetric generating set.
Assuming that $R(\ga)>1/2$ and arguing as in the proof of Proposition~\ref{prop:ue}
conclude that there exists a large-scale Lipschitz map $f\in\LE(\Gamma,\h)$,
an $\epsilon>0$ and an $r_0>0$ such that
\begin{equation}
\label{ass0}
  \rho_f(r) \geq r^{(1+\epsilon)/2}, \quad\text{for all $r\geq r_0$}. 
\end{equation}
Define, for $k\geq 1$, a function $u_k:\Gamma\times\Gamma \to \R$ by
\begin{equation*}
  u_k(s,t) = e^{- \| f(s) - f(t) \|^2 k^{-1}},
         \quad\text{for all $s$, $t\in\Gamma$}.
\end{equation*}
Since the function $\| f(s) - f(t) \|^2$ is of negative type
\cite{delaharpe-valette:T}, each $u_k$ is  positive
definite by Schoenberg's theorem \cite{berg-forst-ptolcag}, and is
also {\it normalized\/} in the sense that $u_k(s,s)=1$, for all
$s\in\ga$.  Further, 
since $f$ is large-scale Lipschitz, the sequence $u_k$ satisfies the
convergence condition \eqref{convergence}.  However, instead of the
support condition (\ref{support}), they possess 
a  weaker decay property. The remainder of the proof will be
devoted to approximating the $u_k$ uniformly by  {\it finite width\/}
positive definite kernels so that both the convergence and support
conditions hold for the approximants.

Recall that the {\it uniform Roe algebra}, $C_u^{*}(\Gamma)$,  is the  
 $C^*$-algebra of bounded operators on $l^2(\Gamma)$ which is the
 norm closure of the subalgebra of
  operators generated by $\op(k)$, where $k$ is a bounded finite width kernel.

\begin{lem}\label{opu}
The operators $\op(u_k)\in C^*_u(\Gamma)$, for all $k\geq 1$.
\end{lem}
\begin{proof}
We show that for every $\kappa>0$ the kernel
$u:\Gamma\times\Gamma\to \C$ defined by
\begin{equation}
\label{eq:u}
  u(s,t) = e^{-\|f(s)-f(t)\|^2\kappa}, 
    \quad\text{$s$, $t\in\ga$}
\end{equation}
defines an element $\op(u)\in C^*_u(\ga)$.  To this end, 
define, for $n\in\N$ 
\begin{equation*}
  k_n(s,t) = \begin{cases} u(s,t), & \text{if $d(s,t)> n$} \\
                         0, & \text{otherwise}. \end{cases}
\end{equation*}
Note that $u-k_n$ is a bounded finite width kernel so that
$\op(u-k_n)\in C^*_u(\ga)$.  Since $\op(u) = \op(u-k_n)+\op(k_n)$ on
compactly supported elements of $l^2(\ga)$, it suffices to show that
$\| \op(k_n) \| \to 0$ as $n\to\infty$.

We proceed using the Schur test.  Since the $k_n$ are non-negative
real-valued and symmetric, it is sufficient to check either one of the
inequalities in (\ref{eqn:schur}).  For this, we will show that there exists
a sequence $C_n\to 0$ such that
\begin{equation*}
\label{eqn:cm}
  \sum_{t\in\Gamma} k_n(s,t) = 
      \sum_{m>n} \sum_{d(s,t)=m} u(s,t) \leq C_n,
     \quad \text{for all $s\in\Gamma$.}
\end{equation*}
This,  in turn, follows from the assertion that there exists $C$
such that
\begin{equation*}
 \sum_{n\geq 0} \sum_{d(s,t)=n} u(s,t) \leq C,
      \quad\text{for all $s\in\Gamma$}.
\end{equation*}
To obtain $C$, let $\sigma$ be the {\it spherical growth function\/} of $\Gamma$
defined by
\begin{equation*}
  \sigma(n) = \card \{\, t\in\Gamma : \text{$d(t,e)=n$} \,\}.
\end{equation*}
Denoting by $S$ the fixed generating set of $\ga$ 
observe that $\sigma(n)\leq (\card S)^n$.  Combining ~\eqref{ass0} and
~\eqref{eq:u} see that if $d(s,t)=n\geq r_0$ then 
$n^{(1+\epsilon)/2} \leq \rho_f(n) \leq \| f(s) - f(t) \|$, and also 
$u(s,t) \leq e^{-\kappa n^{(1+\epsilon)}}$.  Let $m\geq r_0$ be
sufficiently large such that $\card(S) < e^{\kappa m^{\epsilon}}$.
We estimate:
\begin{equation}
\label{bad}
  \begin{aligned}
  \sum_{n\geq 0} \sum_{d(s,t)=n} u(s,t) 
    &= \sum_{n\leq m} \sum_{d(s,t)=n} u(s,t) + 
            \sum_{n>m} \sum_{d(s,t)=n} u(s,t)  \\
    &\leq \sum_{n\leq m} \sigma(n) + 
          \sum_{n>m} \sum_{d(s,t)=n} e^{{-\kappa n^{1+\epsilon}}} \\ 
    &\leq \sum_{n\leq m} \sigma(n) +
          \sum_{n>m}  \sigma(n) e^{{-\kappa n^{1+\epsilon}}} \\
    &\leq \sum_{n\leq m} \sigma(n) +
     \sum_{n>m} \left\{ \frac{\card(S)}{e^{\kappa n^{\epsilon}}} \right\}^n \\
    &\leq \sum_{n\leq m} \sigma(n) +
     \sum_{n>m} \left\{ \frac{\card(S)}{e^{\kappa m^{\epsilon}}} \right\}^n,
  \end{aligned}
\end{equation}
which is both finite and independent of $s\in\Gamma$. We set $C$ equal
to the right hand side of the inequality.  This completes the proof of
the lemma.
\end{proof}

We now complete the proof of the theorem.  Since $u_k$ is normalized
we have $\| \op(u_k) \|\geq 1$.  It is straightforward to show that
since the $u_k$ are positive definite kernels the $\op(u_k)$ are
positive operators.  Let $V_k\in C^*_u(\Gamma)$ be the positive square
root of $\op(u_k)$ and let $W_k\in C^*_u(\Gamma)$ be operators
represented by finite width kernels and such that $\| V_k-W_k\| \| V_k
\| \to 0$.  Define kernels $\widehat u_k$ by
\begin{equation*}
  \widehat u_k (s,t) = \langle\, W_k\delta_t,W_k\delta_s \,\rangle,
     \quad\text{$s$, $t\in\ga$}.
\end{equation*}
The $\widehat u_k$ are positive definite kernels and, since the $W_k$
are represented by finite width kernels, the $\widehat u_k$ are
themselves finite width kernels.  Finally,
\begin{equation*}
  \begin{split}
  | u_k(s,t) - \widehat u_k(s,t) | &= 
    \left| \langle\, \left(\op(u_k)-W_k^*W_k\right)\delta_t,
         \delta_s \,\rangle \right| \\
      &\leq \| V_k^*V_k - W_k^*W_k \| \\
      &\leq \| V_k - W_k \| \left( \| V_k \| + \| W_k \| \right) \\
   &\leq \| V_k - W_k \| \left( 2 \| V_k \| + \| V_k - W_k \| \right),  
  \end{split}
\end{equation*}
which tends to zero as $k\to\infty$.  Consequently $u_k -\widehat u_k\to 0$ 
uniformly on $\Gamma\times \Gamma$ and since the $u_k$
satisfy the convergence condition so do the $\widehat u_k$.
\end{proof}

This theorem has the following interesting consequence.

\begin{thm}
  Let $f: \Gamma \to \h$ be a uniform embedding of a finitely
  generated group into a Hilbert space.  Suppose that $f(\Gamma) \cont
  \h$ is a quasi-geodesic space with the induced metric.  Then
  $C^{*}_{r}(\Gamma)$ is an exact $C^{*}$-algebra.
\end{thm}
\begin{proof}
Since $f$ is a uniform embedding $\Gamma$ is coarsely equivalent to
$f(\Gamma)$ by Proposition~\ref{prop:ue-ce}.  Since both are
quasi-geodesic spaces, it follows from Proposition~\ref{prop:coarse-qi}
that $\Gamma$ is quasi-isometric to $f(\Gamma)$.  But the latter is
isometrically embedded in a Hilbert space, so $R(f(\Gamma)) =1$ by
Proposition~\ref{ex2}.  By quasi-isometry invariance,
Corollary~\ref{cor:qi-inv}, we get $R(\Gamma) = 1$, and hence, by
Theorem~\ref{exact}, $C^{*}_{r}(\Gamma)$ is exact.  
\end{proof}

\begin{rem}
  One might try to argue along the lines of the previous proof to show
  that any uniformly embeddable group is exact.  The first difficulty
  is that if $f(\Gamma)$ is not quasi-geodesic with the induced metric
  then there is no way to deduce that $R(\Gamma) = 1$ even though
  $R(f(\Gamma)) =1$.  However, one need not give up yet.  If one could
  deduce from $R(f(\Gamma)) =1$ the existence of $u_{n}$'s satisfying
  the conditions in Theorem \ref{exact}, then one could pull back the
  $u_{n}$'s to $\Gamma \x \Gamma$ and they would also satisfy the
  necessary condition, thus showing $\Gamma$ was exact.  The
  difficulty here is that the argument used to construct the
  $u_{n}$'s depended on the fact that the spherical growth rate of a
  group with a word length metric is at most exponential--a fact that
  need not hold for an arbitrary discrete metric space like
  $f(\Gamma)$.
\end{rem}

\section{Behavior of $R(\Gamma)$ under direct sums and certain free products}

Let $X$ and $Y$ be metric spaces.  Let  $X\times Y$ be the cartesian
product with the metric
\begin{equation*}
  d_{X\times Y}\left( (x,y),(x',y') \right) = 
       d_X(x,x')+d_Y(y,y').
\end{equation*}
We will obtain a formula for the Hilbert space distortion $R(X \x Y)$
in terms of  $R(X)$ and $R(Y)$.

\begin{prop}
\label{prod}
For metric spaces $X$ and $Y$ we have
$R(X\times Y)=\min \left\{\, R(X),R(Y) \,\right\}$.
\end{prop}
\begin{proof}
Note that,  for fixed $y_0\in Y$, the map $x\longmapsto (x,y_0)$
provides an isometry $X\to X\times Y$.  Applying Theorem~\ref{thm:qii} we 
conclude that $R(X)\geq R(X\times Y)$.  Similarly 
$R(Y)\geq R(X\times Y)$ and so $\min\{\, R(X),R(Y) \,\}\geq R(X\times Y)$.

We must prove the reverse inequality.   Assume that $R(X)\leq R(Y)$. Let
$\eps>0$ be given.  We will show that there exists a large-scale Lipschitz
map $h \cc X\times Y\to\h$ such that $R_h\geq R(X)-\eps$.  From this
one obtains
\begin{align*}
  R(X\times Y) 
    \geq R_h \geq R(X)-\eps 
    = \min \{\, R(X),R(Y) \,\} -\eps,
\end{align*}
and  the desired inequality follows.

According to the definition of $R(X)$ and $R(Y)$ there exist $f \in
\LE(X,\h_X)$ and $g \in\LE(Y,\h_Y)$ such that
\begin{align*}
  R_f &\geq R(X) - \eps \\
  R_g &\geq R(Y) - \eps \geq R(X) - \eps.
\end{align*}
 Define $h\cc X\times Y\to \h=\h_X\oplus \h_Y$ by 
$h(x,y) = f(x)\oplus g(y)$.  From the inequality 
\begin{equation}
\label{inequality2}
  \frac{\alpha+\beta}{\sqrt 2} \leq 
       \left(\alpha^2+\beta^2\right)^{1/2} \leq \alpha+\beta,
     \quad\text{for all $\alpha$, $\beta\geq 0$,}
\end{equation}
we conclude that $h\in \LE(X\times Y,\h)$.
It remains to estimate the compression of $h$, again  using 
(\ref{inequality2}).  We have,
\begin{align*}
  \left\| h(x,y)-h(x',y')\right\| &= 
       \left\| f(x)-f(x') \oplus g(y)-g(y') \right\| \\ 
      &\geq \frac{1}{\sqrt 2} 
 \left\{ \left\| f(x)-f(x')\right\| +\left\| g(y)-g(y') \right\| \right\}
\end{align*}
If $d_{X\times Y}((x,y),(x',y'))\geq r$ then at least one of $d_X(x,x')$
or $d_Y(y,y')\geq r/2$.  Consequently,
\begin{align*}
  \rho_h(r) &= \inf \{\, \|h(x,y)-h(x',y')\| : 
        \text{$d_{X\times Y}((x,y),(x',y'))\geq r$} \,\} \\
     &\geq \frac{1}{\sqrt 2} \inf 
   \left\{\, \left\| f(x)-f(x')\right\| +\left\| g(y)-g(y') \right\| :
        \text{$d_{X\times Y}((x,y),(x',y'))\geq r$} \,\right\} \\
     &\geq \frac{1}{\sqrt 2} \min \left\{\,\rho_f \left(\frac{r}{2}\right),
               \rho_g\left(\frac{r}{2}\right) \,\right\}
\end{align*}
It follows that,
\begin{align*}
  R_h &= \liminf_{r\to\infty} \frac{\log \rho_h(r)}{\log r} \\
      &\geq \liminf_{r\to\infty} \min\left\{\,
         \frac{\log \rho_f\left(\frac{r}{2}\right)}{\log r}, 
       \frac{\log \rho_g\left(\frac{r}{2}\right)}{\log r} \,\right\} \\
      &= \min \{\, R_f,R_g \,\} \geq R(X)-\eps. \qed
\end{align*}
\renewcommand{\qed}{}
\end{proof}

We next study a the free product $\Z * \Z$.  The calculation of $R(\Z * \Z)$
requires a new technique to deform uniform embeddings.  It is
likely that a variant of this technique will apply to other free
products (without amalgam), but we will not not address this in the
present paper. 

\begin{prop}
\label{freegp}
  Let $\F_2$ be the free group on two generators.  Then $R(\F_2)=1$.
\end{prop}
\begin{proof}
Let $X=(V,E)$ be the Cayley graph of $\F_2$, $V\cong\F_2$ being the
set of vertices and $E$ the set of edges.  Let $\h = l^{2}(E)$. Define
\begin{equation*}
  f :\F_2\to \h, \quad
  f(s) = \delta_{e_1(s)} + \cdots + \delta_{e_k(s)},
\end{equation*}
where $\delta_e$ is the Dirac function of the edge $e$ and
$e_1(s),\dots,e_k(s)$ are the edges on the unique path in the Cayley
graph from $s\in\F_2$ to the identity $1\in\F_2$.  Note that
$k=d(s,1)$ so that $\| f(s) \| = \sqrt{d(s,1)}$.  Indeed, the following
assertions can be verified directly:  $\| f(s)-f(t)\|=\sqrt{d(s,t)}$,
for all $s$ and $t\in\ga$ and $\sqrt{r}\leq \rho_f(r)\leq
\sqrt{r+1}$. Hence, the asymptotic compression of $f$ is $1/2$.   

Our strategy for proving the proposition is to produce, by placing
appropriate weights into the above formula for $f$, a family of
large-scale Lipschitz embeddings $f_\eps\in\LE(\F_2,\h)$, for
$0<\eps<1/2$, 
such that $R_{f_\eps}\to 1$ as $\eps\to 1/2$.  Denote
$\xi_\eps(x)=x^\eps$ and define weights by 
$c_{\eps,n} = \xi_\eps(n) = n^\eps$, for $n\in\N$.  Define 
$f_\eps \cc \F_2\to l^2(E)$ by
\begin{equation*}
  f_\eps(s) = c_{\eps,1} \delta_{e_1(s)} +\cdots+ 
         c_{\eps,k} \delta_{e_k(s)},
\end{equation*}
where $k$ and $e_1(s),\dots,e_k(s)$ are as above.  

In order to show that $f_e$ is a large scale Lipschitz map it suffices
to show that there exists $C>0$ such that
\begin{equation*}
  d(s,t) = 1 \Longrightarrow \| f_\eps(s)-f_\eps(t) \|^2 \leq C,
     \quad\text{for all $s$, $t\in\F_2$}.
\end{equation*}
Let $s$, $t\in\F_2$ be such that $d(s,t)=1$.  Denote by $k$ the length
of $s$ and, without loss of generality, $k+1$ the length of $t$.  We
have
\begin{equation*}
  \| f_\eps(s)-f_\eps(t) \|^2 = 
         c_{\eps,1}^2 + (c_{\eps,2}-c_{\eps,1})^2 + \cdots
               +  (c_{\eps,k+1}-c_{\eps,k})^2
\end{equation*}
so that the desired inequality follows from the elementary fact that
$\sum_{j=2}^\infty (c_{\eps,j}-c_{\eps,j-1})^2$ is finite.  Indeed, 
\begin{align*}
  \sum_{j=2}^\infty (c_{\eps,j}-c_{\eps,j-1})^2 &=
    \sum_{j=2}^\infty \left( \int_{j-1}^j \xi'_\eps(x) \,dx \right)^2
 \leq  \sum_{j=2}^\infty \int_{j-1}^j \left(\xi'_\eps(x)\right)^2 \,dx \\
    &= \int_1^\infty \eps^2 x^{2\eps-2} \,dx = \frac{\eps^2}{1-2\eps}.
\end{align*}

To conclude the proof we must show that $R_{f_\eps}\geq 1/2+\eps$.  In
view of Definition~\ref{eqn:rgamma} of the asymptotic compression it
suffices to show that there exists a constant $C_\eps>0$, depending only
on $\eps$, such that 
\begin{equation*}
  \| f_\eps(s) - f_\eps(t) \|^2 \geq C_\eps r^{1+2\eps},
    \quad\text{for all $s$, $t\in\F_2$ with $d(s,t)\geq r$}.
\end{equation*}
Indeed, it follows from this that 
$\rho_{f_\eps}(r) \geq \sqrt{C_\eps} r^{1/2+\eps}$ for all $r\geq 1$ and
hence that $R_{f_\eps}=\liminf_{r\to\infty}
\frac{\log\rho_{f_\eps}(r)}{\log r} \geq 1/2+\eps$.  To prove the
inequality let $s$, $t\in\F_2$ be such that $d(s,t)\geq r$ and assume,
without loss of generality, that $d(1,s)\leq d(1,t)$.  
 Denoting the smallest integer greater than
$r/2$ by $\#(r/2)$, one checks easily that the edges
$e_1(t),\dots,e_{\#(r/2)}(t)$ appear in the expression for $f_\eps(t)$,
but do not appear in that of $f_\eps(s)$.  In particular,
\begin{equation*}
  \| f_\eps(s) - f_\eps(t) \|^2 \geq 
        c^2_{\eps,1} + \cdots + c^2_{\eps,\#(r/2)} \geq
        \int_0^{r/2} \xi_\eps^2(x) \,dx =
         \frac {r^{2\eps+1}}{(2^{2\eps+1})(2\eps+1)}. \qed
\end{equation*}
\renewcommand{\qed}{}\end{proof}

\section{The equivariant case}

Incorporating an action of the group $\Gamma$ the
ideas of the previous section yield results about amenability and
a-T-menability.  To this end, we adapt the previous definitions
and results to the equivariant case.  Let $\Gamma$ be a finitely
generated discrete group, equipped as usual with a word length metric.
Let $X$ be a metric space on which $\Gamma$ acts by isometries.  We
define the equivariant Hilbert space compression of $X$ by restricting
our attention to $\Gamma$-equivariant large-scale Lipschitz maps of
$X$ into Hilbert spaces equipped with actions of $\Gamma$ by {\it
  affine isometries}.  Precisely, define
\begin{equation}
  \LEG(X,\h) = \left\{
  \begin{minipage}{3.0truein}
       $\Gamma$-equivariant large-scale Lipschitz maps 
       $f:X\to \h$, $\h$ a $\Gamma$-Hilbert space;
  \end{minipage} \right.
\end{equation}
the definition of the compression and asymptotic compression of
$f\in\LEG(X,\h)$ are the same as in the non-equivariant case (see
(\ref{eqn:c}) and (\ref{eqn:asc}), respectively); the
$\Gamma$-equivariant Hilbert space compression of $X$ is defined by
\begin{equation*}
  R_\Gamma(X) = \sup\{\, R_{f} : \text{$f \in \LEG(X,\h)$} \,\}.
\end{equation*}
(Compare to Definition~\ref{eqn:rgamma}.)
With these definitions in place the following analogs of
Theorem~\ref{thm:qii} and its corollaries are proved in the same manner.

\begin{thm}
  Let $X$ and $Y$ be metric spaces on which the countable discrete
  group $\ga$ acts by isometries.  If there exists an equivariant
  quasi-isometry $X\to Y$ then $R_\Gamma(X)\geq R_\Gamma(Y)$. \qed
\end{thm}

\begin{cor}
  Let $\ga$ be a finitely generated discrete group.  The invariant
  $R_\Gamma(\Gamma)$ is independent of the finite symmetric generating
  set used to define the length function and metric on $\Gamma$. \qed
\end{cor}
  
Recall that an affine
isometric action of $\ga$ on a Hilbert space $\h$ consists of an
orthogonal representation $t\longmapsto \pi_t$ of $\ga$ on $\h$ and a
function $b:\ga\to\h$ satisfying the {\it cocycle identity \/}
\begin{equation}
\label{eq:cocycle}
  b(st)=\pi_s(b(t))+b(s);
\end{equation}
this identity insures that 
\begin{equation}
\label{eq:aff}
  s\mapsto s\cdot \ :\ga \to \isom(\h),\quad
       s\cdot x = \pi_s(x) + b(s), \quad\text{$x\in\h$},
\end{equation}
defines a
homomorphism from $\ga$ into the group of affine isometries of $\h$.
An affine isometric action is {\it
  metrically proper\/} if for every bounded set $B\subset \h$
the set $\{\; s\in\Gamma : \text{$s\cdot B\cap B \neq \emptyset$} \;\}$
is finite;
equivalently, the cocycle $b$ is {\it proper\/} in the sense that
for every $C>0$ the set $\{\, s\in\ga :\text{$\|b(s)\|\leq C$}\,\}$ is
finite.  A countable discrete group $\Gamma$ has the {\it Haagerup
  property\/} if it admits a metrically proper affine isometric action
on a Hilbert space.  The first part of the next theorem is analogous to
Proposition~\ref{prop:ue}; the second part is analogous to
Theorem~\ref{thm:exact}.

\begin{thm}
\label{hag}
  Let $\Gamma$ be a finitely generated discrete group.  If
  $R_\Gamma(\Gamma)>0$ then $\Gamma$ has the Haagerup property.  If
  $R_\Gamma(\Gamma) > \frac{1}{2}$, then $\Gamma$ is amenable. \qed
\end{thm}

According to the theorem, if a finitely generated
discrete group $\ga$ has an orthogonal representation on a Hilbert
space that admits a cocycle $b$ of sufficiently rapid growth then
it is amenable.  Indeed, suppose that $\pi$ is an orthogonal action
of $\ga$ on $\h$.  A cocycle $b$ for $\pi$ is an element of
$\LEG(\ga,\h)$, where we view $\ga$ as acting on $\h$ by the affine
isometric action (\ref{eq:aff}) and on itself by
multiplication on the left; the required equivariance follows from the
cocycle identity and it is easy to verify that $b$ is large-scale
Lipschitz.  Further, one has
\begin{equation*}
  \| b(s) - b(t) \| = \| \pi_t(b(t^{-1}s) \| = \| b(t^{-1}s) \|
\end{equation*}
from which follows that
\begin{equation*}
  \rho_b(r) = \inf\{\, \| b(s)-b(t) \| : {d(s,t)\geq r} \,\} =
          \inf \{\, \| b(s) \| : {d(s,e)\geq r}\,\}.
\end{equation*}
Thus, if an orthogonal action of $\ga$ on a
Hilbert space $\h$ admits a cocycle $b$ for which 
\begin{equation*}
  R_{b} = \liminf_{r\to\infty}
    \frac{\log \inf \{\, \| b(s) \| : {d(s,e)\geq r} \,\}}{\log r}
   > \frac{1}{2}
\end{equation*}
then it is amenable.  In particular, this is the case if the cocycle
satisfies $\|b(s)\| \geq \left( d(s,e) \right)^{1/2+\eps}$ for
some $\eps>0$.

As an illustration consider once again $\Gamma=\F_2$, the free group
on two generators.  As in the proof of Proposition~\ref{freegp}, let
$X=(V,E)$ be the Cayley graph of $\F_2$, $V\cong\F_2$ being the set of
vertices and $E$ the set of edges.  Let $\h=l^2(E)$
be the Hilbert space of real valued functions
equipped with an orthogonal action, $\pi$, of $\F_2$, and the
function
\begin{equation*}
  b:\F_2 \to l^2(E),\quad
    b(s) = \left\{ 
      \begin{minipage}{2.5truein}
         characteristic function of the set of edges on 
         the unique path from $s$ to the identity
      \end{minipage}\right.
\end{equation*}
satisfies the cocycle identity (\ref{eq:cocycle}).
Consequently, $b\in \LEG(\F_2,l^2(E))$, where we equip $l^2(E)$ with the
affine isometric action (\ref{eq:aff}).

As remarked earlier, 
\begin{equation*}
  \| b(s)-b(t) \| = \sqrt{d(s,t)},
\end{equation*}
for all $s$, $t\in\F_2$, and the asymptotic compression of $b$ is 
$\displaystyle{ R_b = \liminf_{r \to \infty} 
   \frac{\log\rho_{b}(r)}{\log r} } = 1/2$.  In particular,  
the equivariant Hilbert space compression of $\F_2$ satisfies
$R_{\F_2}(\F_2)\geq 1/2$.  On the other hand, since $\F_2$ is not
amenable we have $R_{\F_2}(\F_2)\leq 1/2$.  Hence $R_{\F_2}(\F_2)
=1/2$.   This should be compared to Proposition~\ref{freegp}, in which
we proved, by deforming the cocycle $b$, that that $R(\F_2)=1$. 

\section{Appendix}

In this appendix we review several known relations between uniform
embeddings and other notions of coarse geometry \cite{ roe:cbms}.  We
include some of the elementary proofs for the convenience of the reader.

Let $X$ and $Y$ be metric spaces.  A {\it coarse map\/} is a function
$f: X \to Y$ satisfying the following two conditions: 
\begin{ilist}
\item  For every $R>0$ there exists an $S > 0$ such that
  \begin{equation*}
    d_{X}(x,x') \leq R \Longrightarrow d_{Y}(f(x), f(x')) \leq S.
  \end{equation*}
\item  If $B \cont Y$ is bounded, then $f^{-1}(B)$ is bounded. 
\end{ilist}
A coarse map $f: X \to Y$ is a {\it coarse equivalence\/} if there is a
coarse map $g:Y \to X$ and a $K > 0$ such that
\begin{equation}
\begin{aligned}
\label{eq:coarseK}
  d(gf(x),x) &\leq K,\quad\text{for all $x\in X$}, \\
  d(fg(y),y) &\leq K,\quad\text{for all $y\in Y$.}    
\end{aligned}
\end{equation}

\begin{prop}
Let $X$ and $Y$ be metric spaces. A function $f: X \to Y$ is a uniform
embedding if and only if it satisfies the following two conditions:
\begin{ilist}
\item  For every $R>0$ there exists an $S > 0$ such that
  \begin{equation*}
      d_{X}(x,x') \leq R \Longrightarrow d_{Y}(f(x), f(x')) \leq S.
  \end{equation*}
\item  For every $S>0$ there exists an $R > 0$ such that
  \begin{equation*}
      d_{X}(x,x') \geq R \Longrightarrow d_{Y}(f(x), f(x')) \geq S.
  \end{equation*}
\end{ilist} 
\end{prop}

\noindent
Condition (ii) in this proposition implies condition (ii) in the
preceding definition so that a uniform embedding is a coarse map.
Conversely, coarse map need not be a uniform embedding.

\begin{proof}[Sketch]
Let $f:X\to Y$ be a function satisfying (i) and (ii).  By virtue of
(i) we may define a non-decreasing, real-valued function $\rho_+$ by
  \begin{equation}
    \label{eq:rho+}
    \rho_+(r) = \sup_{d_X(x,y)\leq r} d_Y(f(x),f(y));
  \end{equation}
we define $\rho_-=\rho_f$ according to (\ref{eqn:c}).  These functions
satisfy the inequalities on (\ref{eqn:ue}).  By virtue of (ii), $\rho_-$
is proper, as is $\rho_+$.

We omit verification of the converse.
\end{proof}

\begin{prop}
\label{prop:ue-ce}
  Let $X$ and $Y$ be metric spaces. A function $f: X
  \to Y$ is a uniform embedding if and only if it is a coarse
  equivalence of $X$ with $f(X) \cont Y$ with the induced metric. \qed
\end{prop}
   
\begin{prop}
\label{prop:coarse-qi}
  Let $X$ and $Y$ be  {quasi-geodesic} metric spaces.  A
  function $f: X \to Y$ is a coarse equivalence if and only if it is a
  quasi-isometric equivalence.  
\end{prop}
\begin{proof}
A quasi-isometric equivalence is always a coarse equivalence; we show
the converse under the assumption that $X$ and $Y$ are quasi-geodesic
spaces.

Let $f:X\to Y$ be a coarse equivalence.  The function $\rho_+$ defined
in (\ref{eq:rho+}) satisfies the second inequality in (\ref{eqn:ue});
hence, according to Proposition~\ref{qg}, $f$ is large-scale
Lipschitz.  It remains only to find constants satisfying the
first inequality in (\ref{eq:qidef}).

Let $g:Y\to X$ be a coarse map
satisfying (\ref{eq:coarseK}) for some $K\geq 0$.  Arguing as for $f$
conclude that $g$ is large-scale Lipschitz.  Let $C>0$ and
$D\geq 0$ be constants as in the definition (\ref{eq:lslip}) of
large-scale Lipschitz; that is, satisfying 
\begin{equation*}
  d_X(g(y),g(y'))\leq C d_Y(y,y') +D, 
\end{equation*}
for all $y$, $y'\in Y$.
Let $x$ and $x'\in X$ and calculate
\begin{equation*}
\begin{aligned}
    d_X(x,x') &\leq d_X(x,gf(x)) + d_X(gf(x),gf(x'))+d_X(gf(x'),x') \\
              &\leq 2K + C d_Y(f(x),f(x')) +D
\end{aligned}
\end{equation*}
from which follows that
\begin{equation*}
  C^{-1} d_X(x,x') - C^{-1}(2K+D) \leq d_Y(f(x),f(x')).  \quad\qed
\end{equation*}
\renewcommand{\qed}{}
\end{proof}

\begin{prop}
  Let $X$ and $Y$ be discrete metric spaces with $X$ quasi-geodesic.
  A uniform embedding $f:X \to Y$ is a quasi-isometry if and only if
  $f(X)$, with the metric induced from $Y$, is a quasi-geodesic space.
  \qed
\end{prop}
\begin{proof}
Let $f:X\to Y$ be a uniform embedding and assume that $f(X)$ is a
quasi-geodesic space.  By Proposition~\ref{prop:ue-ce}, $f$ is a
coarse equivalence between $X$ and $f(X)$, and by
Proposition~\ref{prop:coarse-qi} it is  
a quasi-isometry.

For the converse, let $f:X\to Y$ be a quasi-isometry and   let $C>0$ and
$D\geq 0$ be associated constants.   Let
$\delta>0$ and $\lambda\geq 1$ be the constants as in 
(\ref{eqn:qg}) reflecting the fact that $X$ is a quasi-geodesic space.
We show that  two points in $f(X)$ are connected by a
sequence of points satisfying  (\ref{eqn:qg}) with $\delta$
and $\lambda$ replaced by
\begin{equation*}
  \delta' = \dfrac{3\delta C}{2} +D, \qquad
  \lambda' = \dfrac{2 C(D+\delta')\lambda}{\delta} +1.
\end{equation*}

Let $f(x)$ and $f(y)\in f(X)$.  We may assume that $d_Y(f(x),f(y))\geq
\delta'$, for if not $f(x),f(y)$ is the required sequence.  As in the
proof of Proposition~\ref{qg} we obtain a sub-sequence of points in $X$,
$x_{0},x_{1},\dots,x_{m},y$ connecting $x$ and $y$ and satisfying the
required estimates. (The subsequence indices are suppressed here for
easier reading.)  One can show directly that
\begin{equation}
\label{eq:simple}
  d_Y(f(x_{j-1}),f(x_{j})) \leq \delta', \quad
     d_Y(f(x_{m}),y)\leq \delta',
\end{equation}
so it remains  to verify that
\begin{equation}
\label{eq:finalsum}
  \sum_{j=1}^m d_Y(f(x_{j-1}),f(x_{j})) + d_Y(f(x_{m}),f(y))
          \leq \lambda' d_Y(f(x),f(y)).
\end{equation}
From the first inequality in (\ref{eq:qidef}) we conclude that
$d_X(x,y)\leq C (d_Y(f(x),f(y))+D)$.  Since
$d_Y(f(x),f(y))\geq\delta'$ one concludes that
\begin{equation*}
  d_X(x,y) \leq C\left(1+\dfrac{D}{\delta'}\right) d_Y(f(x),f(y)).
\end{equation*}
Combining this inequality with (\ref{eq:simple}) and the bound on $m$
from  \eqref{qg}, the sum in
(\ref{eq:finalsum}) is bounded by  
\begin{equation*}
  (m+1)\delta' \leq 
  \left( \dfrac{2\lambda}{\delta}\, d_X(x,y) + 1 \right)\delta'
             \leq \dfrac{2\lambda}{\delta}\, 
   C\left( \delta'+ {D} \right) d_Y(f(x),f(y)) + \delta' \leq \lambda' d_Y(f(x),f(y)),
\end{equation*}
where we  again use the assumption that $d_Y(f(x),f(y))\geq\delta'$.  This concludes the proof.
\end{proof}

The inclusion of a finitely generated group as a subgroup in another
finitely generated group is a uniform embedding, but it's range need not
be a quasi-geodesic metric space with the induced metric; the inclusion
of $\Z$ in the discrete 3-dimensional Heisenberg group provides an
example of this phenomenon.

\end{document}